\documentclass[11pt]{article}

\usepackage{amsfonts,amsmath,amssymb,amsthm}
\usepackage[all]{xypic}

\newcommand{\Cc}{\mathbb{C}}  

\newcommand{\Zz}{\mathbb{Z}}
\renewcommand{\epsilon}{\varepsilon}

\renewcommand{\ge}{\geqslant}

\renewcommand {\geq}{\geqslant}

\newcommand{\Sing}{\mathop{\mathrm{Sing}}\nolimits}
\newcommand{\rk}{\mathop{\mathrm{rk}}\nolimits}

\newcommand{\grad}{\mathop{\mathrm{grad}}\nolimits}

\newcommand{\Baff}{{\mathcal{B}_\mathit{\!aff}}}
\newcommand{\Binf}{{\mathcal{B}_\infty}}
\newcommand{\B}{{\mathcal{B}}}
\newcommand{\Fgen}{F_\mathit{\!gen}}
\newcommand{\Fgenb}{\bar{F}_\mathit{\!gen}}
\newcommand{\Frond}{F^\circ}

\newcommand{\Ddic}{D_{\!\textit{dic}}}

\newcommand{\Dic} {\mathrm{Dic}}

\newcommand{\m}{\mathfrak{m}}

\newcommand{\Faff}{{F_{\!\mathit{aff}}}}
\newcommand{\Taff}{{T_{\!\mathit{aff}}}}
\newcommand{\jaff}{{j_{\mathit{aff}}}}
\newcommand{\jrond}{j^\circ}

\newcommand{\Ker}{\mathop{\mathrm{Ker}}\nolimits}
\newcommand{\id}{\mathop{\mathrm{id}}\nolimits}

\renewcommand{\Im}{\mathop{\mathrm{Im}}\nolimits}

\newtheorem*{theorem*}{Theorem}  
\newtheorem{theorem}{Theorem}    
\newtheorem{lemma}[theorem]{Lemma}        
\newtheorem{proposition}[theorem]{Proposition}
\newtheorem*{theoremloc*}{Local invariant cycles theorem}
\newtheorem*{theoremglob*}{Global invariant cycles theorem}
\theoremstyle{remark}
\newtheorem*{remark*}{Remark}

\begin{document}

\title{\textbf{Irregular fibers of complex polynomials in two variables}}
\author{\textbf{\textsf{Arnaud Bodin}}}
\date{}
\maketitle

\section*{Introduction}
\label{sec:intro}

Let $f : \Cc^n \longrightarrow \Cc$ be a polynomial. The 
\emph{bifurcation set} $\B$ for $f$ is the minimal set of points of $\Cc$
such that $f : \Cc^n \setminus f^{-1}(\B) \longrightarrow \Cc \setminus \B$ is
a locally trivial fibration. 
For $c\in \Cc$, we denote the fiber $f^{-1}(c)$ by $F_c$.
The fiber $F_c$ is \emph{irregular} if $c$ is in $\B$.
If  $s\notin \B$, then $F_s$  is a \emph{generic fiber}  and is
denoted by $\Fgen$. 
The \emph{tube} $T_c$ for the value $c$ is a neighborhood $f^{-1}(D^2_\epsilon(c))$ 
of the fiber 
$F_c$, where $D^2_\epsilon(c)$ stands for a $2$-disk in $\Cc$, centered at $c$, of radius 
$\epsilon \ll 1$.
We assume that affine critical singularities are isolated.
The value $c$ is \emph{regular at infinity} if there exists a compact set $K$ of $\Cc^2$ such
that the restriction of $f$ to $T_c \setminus K \longrightarrow D^2_\epsilon(c)$ is
a locally trivial fibration.
Set $n=2$.
Let $j_c : H_1(F_c) \longrightarrow H_1(T_c)$ be the morphism induced by the inclusion
of $F_c$ in $T_c$. One of the consequences of our study of $j_c$ is 
the following :

\begin{theorem*}
$j_c$ is an isomorphism if and only if $c$ is a regular value at infinity.
\end{theorem*}

E.~Artal-Bartolo, Pi.~Cassou-Nogu\`es and A.~Dimca has proved this result in \cite{ACD}
for polynomials with a connected fiber $F_c$. But for example the fiber
$F_0$ of Broughton's polynomial $f(x,y) = x(xy+1)$ is not connected.
We also give necessary and sufficient conditions for $j_c$ to be injective and
surjective: 
let $G_c$ the dual graph of $F_c = f^{-1}(c)$ and $\bar G_c$ the dual graph 
of a compactification of the fiber $F_c$ obtained by a 
resolution at infinity of $f$. We find that $\rk \Ker j_c =
n(F_c)-1 + \rk H_1(\bar G_c)-\rk H_1(G_c)$ where $n(F_c)$ is the number of 
connected components of $F_c$. So  $j_c$ is injective if and only if $F_c$ is connected and 
$H_1(G_c)$ is isomorphic to $H_1(\bar G_c)$.

\bigskip

We apply these results to the study of neighborhoods of irregular fibers.
Set $n \geq 2$.
Let $\Frond_c$ be the \emph{smooth part} of $F_c$: 
$\Frond_c$ is obtained by intersecting $F_c$
with  a large $2n$-ball and cutting out a small neighborhood
of the (isolated) singularities.
Then $\Frond_c$ can be embedded in $\Fgen$.
We study the following commutative diagram that links the three
elements $\Frond_c$, $\Fgen$, and $T_c$: 
$$
\xymatrix{
{}H_q(\Frond_c)\ar[r]^{j_c^\circ}\ar[d]_{\ell_c}    &H_q(T_c)  \\
  H_q(\Fgen)    \ar[ur]_{k_c}   & {}
}
$$
where $\ell_c$ is  the morphism
induced in integral homology by the embedding;  
 $j^\circ_c$ and  $k_c$ are induced by inclusions.
The morphism $k_c$ is well-known and
$V_q(c) = \Ker k_c$ are \emph{vanishing cycles} for the value $c$.
Let $h_c$ be the monodromy induced on $H_q(\Fgen)$ by a small
circle around the value $c$.
Then we prove that the image of $\ell_c$ are invariant cycles by $h_c$: 
$$\Ker (h_c - \id) = \ell_c \big( H_q(\Frond_c) \big).$$
This formula for the case $n=2$ has been obtained by 
F.~Michel and C.~Weber in \cite{MW}.

Finally we give a description of vanishing cycles with respect
to eigenvalues of $h_c$ for homology with complex coefficients. 
For $\lambda \not= 1$ and $p$ a large integer the characteristic
space $E_\lambda = \Ker (h_c - \lambda \id)^p$ is composed of vanishing
cycles for the value $c$.
For $\lambda = 1$ the situation is different.
If $K_q(c) = V_q(c) \cap \Ker(h_c-\id)$ are 
invariant and vanishing cycles we have 
$$K_q(c) = \ell_c \big(\Ker j^\circ_c\big).$$
And for  $n=2$ we get the formula
$$\rk K_1(c) = r(F_c)-1 + \rk H_1(\bar{G}_c).$$

\section{Irregular fibers and tubes}
\label{sec:tub}

\subsection{Bifurcation set}

We can describe the bifurcation set $\B$ as follows: let 
$\Sing = \left\{ z \in \Cc^n\ | \ \grad_f(z)=0 \right\}$
be the set of \emph{affine critical points}
and let $\Baff = f(\Sing)$
be the set of \emph{affine critical values}.
The set $\Baff$ is a subset of $\B$. 
The value $c\in \Cc$ is \emph{regular at infinity} if there
exists a disk ${D}$ centered at $c$ and a compact set $K$
of $\Cc^n$
with a locally trivial fibration 
$f : f^{-1} \left( {D} \right)\setminus {K}
\longrightarrow {D}$. The non-regular
values at infinity are the \emph{critical values at infinity} and are collected
in $\Binf$. The finite set $\B$ of critical values is now:
$$\B = \Baff \cup \Binf.$$
In this article we always assume that \textbf{affine singularities are isolated}, that is to say
that $\Sing$ is an isolated set in $\Cc^n$.
For $n=2$ this hypothesis implies that the generic fiber is a connected set.

\subsection{Statement of the result}
In this paragraph $n=2$.
The inclusion of $F_c$ in $T_c$ induces a morphism $j_c : H_1(F_c) \longrightarrow H_1(T_c)$.
\begin{theorem}
\label{th:jnonbij}
 $j_c : H_1(F_c) \longrightarrow H_1(T_c)$ is an isomorphism if and only if 
 $c\notin \B_\infty$.
\end{theorem}

When $F_c$ is a connected fiber this result has been obtained in \cite{ACD}.
We generalize this study and we give simple criteria to determine
whether $j_c$ is injective or surjective.
We firstly recall notations and results from \cite{ACD}.

\bigskip

Let denote $\Faff =  {F_c \cap B^4_R}$ ($R\gg 1$) and $F_\infty = \overline{F_c \setminus \Faff}$,
thus $\Faff \cap F_\infty = K_c = f^{-1}(c) \cap S^3_R$ is the \emph{link at infinity}
for the value $c$.
Similarly $\Taff =  {T_c \cap B^4_R}$ and $T_\infty = \overline{T_c \setminus \Taff}$.
We denote $j_\infty : H_1(F_\infty) \longrightarrow H_1(T_\infty)$ the morphism
induced by inclusion.
The morphism $\jaff : H_1(\Faff) \longrightarrow H_1(\Taff)$ is an isomorphism.
 $H_1(\Faff\cap F_\infty)$ and 
 $H_1(\Taff\cap T_\infty)$ are isomorphic.

\bigskip

Mayer-Vietoris exact sequences for the decompositions
$F_c=\Faff \cup F_\infty$ and $T_c= \Taff \cup T_\infty$ give
the commutative diagram $(\mathcal{D})$:
\begin{equation*}
\label{diag:mv}
\xymatrix{
 0\ar[r] & H_1(F_\infty\cap\Faff) \ar[r]^-g\ar[d]^\cong & H_1(F_\infty)\oplus H_1(\Faff) 
    \ar[r]^-h\ar[d]^{j_\infty \oplus \jaff} & H_1(F_c) \ar[r]\ar[d]^{j_c} & 0 \\
 0\ar[r] & H_1(T_\infty\cap\Taff) \ar[r]^-{g'} & H_1(T_\infty)\oplus H_1(\Taff) \ar[r]^-{h'} & 
        H_1(T_c) \ar[r] & H_0(T_\infty\cap\Taff). }
\end{equation*}

The $0$ at the upper-right corner is provided by the injectivity of
 $H_0(F_\infty \cap \Faff) \longrightarrow H_0(F_\infty)$ 
($F_c$ need not to be a connected set) hence
 $H_0(F_\infty \cap \Faff) \longrightarrow
H_0(F_\infty)\oplus H_0(\Faff)$ is injective.

\subsection{Resolution of singularities}

To compactify the situation, for $n=2$, we need resolution of singularities at infinity
\cite{LW}:
$$
\xymatrix{
{}\Cc^2   \ar[r] \ar[d]_-f      &\Cc P^2 \ar[d]_-{\tilde{f}}    
    &\Sigma_w \ar[l]_-{\pi_w} \ar[dl]^-{\phi_w} \\
\Cc \ar[r]      &\Cc P^1
}
$$
 $\tilde{f}$ is the map
coming from the homogenization of $f$; $\pi$ is the minimal blow-up
of some points on the line at infinity $L_\infty$ of $\Cc P^2$
in order to obtain a well-defined morphism $\phi_w : \Sigma_w \longrightarrow \Cc P^1$:
this the \emph{weak resolution}. 
We denote $\phi_w^{-1}(\infty)$ by  $D_\infty$,
and let $\Ddic$ be the set of components $D$ of $\pi_w^{-1}(L_\infty)$
that verify $\phi_w(D) = \Cc P^1$. Such a $D$ is a \emph{dicritical component}.
The \emph{degree} of a dicritical component $D$ is the degree of the branched covering
$\phi_w : D \longrightarrow \Cc P^1$.
For the weak resolution the divisor $\phi_w^{-1}(c) \cap \pi^{-1}(L_\infty)$, $c \in \Cc$,
is a union of {bamboos} (possibly empty) (a \emph{bamboo} is a divisor whose dual graph is a linear tree).
The set $\B_\infty$ is the set of values of $\phi_w$ on non-empty bamboos with the set of critical 
values of the restriction of $\phi_w$ to the dicritical components.

We can blow-up more points to obtain the \emph{total resolution},
$\phi_t : \Sigma_t \longrightarrow \Cc P^1$, such that all fibers  of 
$\phi_t$ are normal crossing divisors that intersect the dicritical components transversally;
moreover we blow-up affine singularities.
Then $D_\infty = \phi_t^{-1}(\infty)$ is the same as above and for
$c \in \B$ we denote $D_c$ the divisor $\phi_t^{-1}(c)$.

The \emph{dual graph} $\bar{G}_c$ of $D_c$ is obtained as follows:
one vertex for each irreducible component of $D_c$ and one edge between two
vertices for one intersection of the corresponding components. 
A similar construction is done for $D_\infty$, we know
that $\bar G_\infty$ is a tree \cite{LW}.
The \emph{multiplicity} of a component is the multiplicity of 
$\phi_t$ on this component.

\subsection{Study of $j_\infty$}

See \cite{ACD}.
Let $\phi$ be the weak resolution map for $f$. Let denote by $\Dic_c$ the set of points 
$P$ in the dicritical components, such that $\phi(P)=c$.
To each $P\in \Dic_c$ is associated one, and only one, connected component $T_P$ of $T_\infty$;
$T_P$ is the \emph{place at infinity} for $P$.
We have $T_\infty = \coprod_{P\in \Dic_c} T_P$ and we set
 $F_P = T_P \cap F_\infty = T_P \cap F_c$ and $K_P = \partial F_P$, 
finally $n(F_P)$  denotes the number of connected components of $F_P$. Let $\bar{F}_P$
be the strict transform of $c$ by  $\phi$, intersected with $T_P$.
The study of $j_\infty$ follows from the study of $j_P : H_1(F_P) \longrightarrow H_1(T_P)$. 
Let $\m_P$ be the intersection multiplicity of $\bar{F}_P$ with the divisor 
$\pi_w^{-1}(L_\infty)$ at $P$.

\paragraph{Case of  $P \in \bar{F}_P$.} 
The group $H_1(T_P)$ is isomorphic to $\Zz$ and is generated by $[M_P]$,
 $M_P$ being the boundary of a small disk with transversal intersection 
with the dicritical component.
Moreover if $F_P = \coprod_{i=1}^{n(F_P)} {F_P^i}$ then $j_P([F_P^i]) = j_P([K_P^i]) = \m_P^i[M_P]$.

\paragraph{Case of $P$ being in a bamboo.}
The group  $H_1(T_P)$ is also isomorphic to $\Zz$   and is generated 
by $[M_P]$, $M_P$ being the boundary of a small disk, with transversal intersection
with the last component of the bamboo.
Then  $j_P[F_P^i] = j_P[K_P^i] = \m_P^i.\ell_i[M_P]$.
The integer $\ell_i$ only depends  of the position where  $F_P^i$ intersects the bamboo,
moreover $\ell_i \ge 1$ and $\ell_i = 1$ if and only if $F_P^i$ intersects the bamboo 
at the last component.
For a computation of $\ell_i$, refer to \cite{ACD}.

\bigskip

As a consequence $j_P$ is injective if and only if $n(F_P) = 1$ and 
$j_\infty$ is injective if and only if $n(F_P)=1$ for all $P$ in $\Dic_c$.
In fact the rank of the kernel of $j_\infty$ is the sum of the ranks of the kernels of  $j_P$
then
$$\rk \ker j_\infty = \sum_{P\in \Dic_c} (n(F_P)-1).$$

Finally $j_\infty$ is surjective if and only if for all $P\in\Dic_c$, $j_P$ is 
surjective.

\subsection{Acyclicity}

The value $c$ is \emph{acyclic} if the morphism
$\psi : H_0(T_\infty\cap\Taff)\longrightarrow H_0(T_\infty)\oplus H_0(\Taff)$ 
given by the Mayer-Vietoris exact sequence is injective.

\bigskip

Let give some interpretations of the acyclicity condition.
\begin{enumerate}
  \item The injectivity of $\psi$ can be view as follows:
two branches at infinity that intersect the same place at infinity have to be in 
different connected components of $F_c$.
   \item Let $G_c$ be the  \emph{dual graph} of $F_c$ 
(one vertex for an irreducible component of $F_c$, two vertices are joined by an edge 
if the corresponding irreducible components have non-empty intersection, 
if a component has auto-intersection it provides a loop)
and let $G_{c,P}$ be the graph obtained from $G_c$ by 
adding edges to vertices that correspond to
the same place at infinity $T_P$. 
In other words $c$ is acyclic if and only if there is no new cycles in
  $G_{c,P}$, that is to say $H_1(G_c) \cong  H_1(G_{c,P})$ for all $P$ in $\Dic_c$.
   \item Another interpretation is the following:  $c$  is acyclic if and only if the morphism $h'$ 
of the diagram $(\mathcal{D})$ is surjective. This can be proved by the exact sequence:
\begin{align*}
\ H_1(T_\infty)\oplus& H_1(\Taff)\stackrel{h'}{\longrightarrow}
H_1(T_c)\stackrel{\varphi}{\longrightarrow} 
 H_0(T_\infty\cap\Taff) \stackrel{\psi}{\longrightarrow} \\
& \stackrel{\psi}{\longrightarrow}
 H_0(T_\infty)\oplus H_0(\Taff) \longrightarrow H_0(T_c).
\end{align*}
   \item Let consider the above Mayer-Vietoris exact sequence in reduced homology, the  morphism
$\widetilde{\psi} : \widetilde{H}_0(T_\infty\cap\Taff) \longrightarrow 
\widetilde{H}_0(T_\infty)\oplus \widetilde{H}_0(\Taff)$ is surjective because
$\widetilde{H}_0(T_c) = \{ 0\}$. Moreover $\widetilde{\psi}$ is injective 
if and only if  $\psi$ is injective.
As $\widetilde{\psi}$ is  surjective,  $\widetilde{\psi}$  is injective if and only if
 $\rk \widetilde{H}_0(T_\infty\cap\Taff)=
 \rk \widetilde{H}_0(T_\infty) +\rk \widetilde{H}_0(\Taff) $,
that is to say $c$ is acyclic if and only if 
\begin{equation}
\label{eq:dim}
\sum_{P\in \Dic_c}{\!\! n(F_P)} \ \  -\  1 =  \# \Dic_c  -1 + n(F_c)-1.
\tag{$\star$}
\end{equation}
\end{enumerate}

This implies the lemma:
\begin{lemma}
\label{lem:acy}
$j_\infty \text{ is injective } \Longleftrightarrow  F_c  \text{ is a connected set and }
c \text{ is acyclic.}$
\end{lemma}

\begin{proof}
If $j_\infty$ is injective then  $n(F_P)= 1$ for all  $P$ in $\Dic_c$, 
then $H_0(T_\infty\cap\Taff) \cong H_0(T_\infty)$ and  $\psi$ is injective, hence 
$c$ is acyclic and from equality (\ref{eq:dim}), we have
$n(F_c)=1$ \emph{i.e.} $F_c$ is a connected set.
Conversely, if $c$ is acyclic and $n(F_c)=1$ then equality (\ref{eq:dim}) gives
$n(F_P)=1$ for all $P$ in $\Dic_c$, thus $j_\infty$ is injective.
\end{proof}

\bigskip
Let us define a stronger notion of acyclicity.
Let $\bar{G}_c$ be the dual graph of $\phi^{-1}(c)$.
The graph  $\bar{G}_c$ can be obtained from $G_c$ by adding  edges between vertices
that belong to the same place at infinity
for all $P$ in $\Dic_c$.
The value $c$ is \emph{strongly acyclic} if  $H_1(\bar{G}_c) \cong  H_1(G_{c})$.
Strong acyclicity implies acyclicity, but the converse can be false.
However if $F_c$ is a connected set (that is to say $G_c$ is a connected graph)
then both conditions are equivalent. 
This is implicitly expressed in the next lemma, which is 
just a result involving graphs.
\begin{lemma}
\label{lem:graph}
$\displaystyle 
\rk H_1(\bar{G}_c) - \rk H_1(G_c) 
  =  \sum_{P\in \Dic_c} \big(n(F_P)-1\big) \ -\  \big(n(F_c)-1\big).$
\end{lemma}

\subsection{Surjectivity}
\label{subsec:surj}

\begin{proposition}
\label{prop:jsurj}
$j_c \text{ surjective  } \Longleftrightarrow 
  j_\infty \text{ surjective and } c \text{ acyclic.}$
\end{proposition}

\begin{proof}
Let us suppose that $j_c$ is surjective then a version of the five lemma applied to diagram 
$(\mathcal{D})$ proves that
$j_\infty$ is surjective. As $j_c$ and $j_\infty$ are surjective, diagram
$(\mathcal{D})$ implies that $h' : H_1(T_\infty) \oplus H_1(\Taff) \longrightarrow H_1(T_c)$ 
is surjective, that means that $c$ is acyclic.

Conversely if $j_\infty$ is surjective and $c$ is acyclic then $h'$ is
surjective and diagram $(\mathcal{D})$ implies that $j_c$ is surjective.
\end{proof}

\subsection{Injectivity}
\label{subsec:inj}

\begin{proposition}
$j_c \text{ is injective } \Longleftrightarrow 
  F_c \text{ is a connected set and } c \text{ is acyclic.}$
\end{proposition}
It follows from lemma \ref{lem:acy} and from the next lemma.
\begin{lemma} 
\label{lem:jinj}
$j_c \text{ injective } \Longleftrightarrow j_\infty \text{ injective. }$\\
Moreover the rank of the kernel is:
\begin{align*}
 \rk \ker j_c &= \rk \ker j_\infty =  \sum_{P\in \Dic_c} \big(n(F_P)-1\big) \\
             &=n(F_c)-1 + \rk H_1(\bar{G}_c) - \rk H_1(G_c).
\end{align*}
\end{lemma}

\begin{proof}
The first part of this lemma can be proved by a version of the five
lemma. However we shall only prove the equality of the ranks of
$\ker j_c$ and $\ker j_\infty$. 
It will imply the lemma because we already know that 
$\rk \ker j_\infty = \sum_{P\in \Dic_c} \big(n(F_P)-1\big)$
and from lemma  \ref{lem:graph} we then have 
$\rk \ker j_\infty =n(F_c)-1 + \rk H_1(\bar{G}_c) - \rk H_1(G_c)$.

\bigskip

The study of the morphism $j_c : H_1(F_c)\longrightarrow H_1(T_c)$
is equivalent to the study of the morphism $H_1(\Taff) \longrightarrow H_1(T_c)$
induced by inclusion that, by abuse, will also be denoted by $j_c$.
To see this, it suffices to remark that $F_c$ is obtained from $\Faff = F_c\cap B^4_R$ 
by gluing $F_c\cap S^3_R \times [0,+\infty[$ to its boundary $F_c\cap S^3_R $.
Then the morphism $H_1(\Faff) \longrightarrow H_1(F_c)$ induced by inclusion is an isomorphism;
finally $\jaff : H_1(\Faff) \longrightarrow H_1(\Taff)$ is also an isomorphism.
The long exact sequence for the pair $(T_c,\Taff)$ is:
$$ H_2(T_c) 
\longrightarrow H_2(T_c,\Taff) \longrightarrow H_1(\Taff)
 \stackrel{j_c}{\longrightarrow} H_1(T_c)
$$
but $H_2(T_c) = 0$ (see \cite{ACD} for example)
then the rank of $\ker j_c$ is the rank of $ H_2(T_c,\Taff)$.

On the other hand, the study of  $j_\infty : H_1(F_\infty) \longrightarrow H_1(T_\infty)$
is the same as the study of $H_1(\partial T_\infty) \longrightarrow H_1(T_\infty)$ 
induced by inclusion (and denoted by $j_\infty$) because the morphisms 
$H_1(\partial F_\infty) \longrightarrow H_1(F_\infty)$ and
$H_1(\partial F_\infty) \longrightarrow H_1(\partial T_\infty) $  induced by inclusions
are isomorphisms. The long exact sequence for $(T_\infty, \partial T_\infty)$ is:
$$ H_2(T_\infty) 
\longrightarrow H_2(T_\infty,\partial T_\infty) \longrightarrow H_1(\partial T_\infty)
 \stackrel{j_\infty}{\longrightarrow} H_1(T_\infty)
.$$
As $H_2(T_\infty) = 0$ (see \cite{ACD}), then the rank of $\ker j_\infty$ is the same as 
$H_2(T_\infty,\partial T_\infty)$.

Finally the groups $H_2(T_\infty,\partial T_\infty)$ and
$H_2(T_c,\Taff)$ are isomorphic by excision, and then the ranks of $\ker j_c$ and 
of  $\ker j_\infty$ are equal. That completes the proof.
\end{proof}

\subsection{Proof of the theorem}

If $c\notin \B_\infty$, then the isomorphism
 $\jaff : H_1(\Faff) \longrightarrow H_1(\Taff)$
implies that $j_c$ is an isomorphism.
Let suppose that $c$ is a critical value at infinity and that $j_c$ is injective.
We have to prove that $j_c$ is not surjective.
As $j_c$ is injective then by lemma  \ref{lem:jinj},
$j_\infty$ is injective. 
By proposition \ref{prop:jsurj} it suffices to prove that $j_\infty$ is not 
surjective.
Let $P$ be a point of $\Dic_c$  that provides irregularity at infinity
for the value $c$, then  $n(F_P) = 1$ because $j_\infty$
is injective. 
Let us prove that the morphism $j_P$ is not surjective.
For the case of $P \in \bar{F}_P$, the intersection multiplicity $\m_P$ is greater than $1$,
then $j_P$ is not surjective.
For the second case, in which $P$  belongs to a bamboo, then
$\m_P.\ell_i > 1$ except for the situation where only one strict transform intersects
the bamboo at the last component. This is exactly the situation excluded 
by the lemma ``bamboo extremity fiber'' of \cite{MW}.
Hence $j_\infty$ is not surjective and $j_c$ is not an isomorphism.

\subsection{Examples}
 
We apply the results to two classical examples.

\paragraph{Broughton polynomial.} Let $f(x,y) =x(xy+1)$, then 
$\Baff = \varnothing$, $\B = \Binf = \{0\}$. Then for $c\not=0$,
$j_c$ is an isomorphism. The value $0$ is acyclic since 
$H_1(G_0) \cong H_1(\bar G_0)$.
The fiber $F_0$ is not connected hence $j_0$ is not injective.
As the new component of $\bar G_0$ is of multiplicity $1$ the corresponding
morphism $j_\infty$ is surjective, hence $j_0$ is surjective.
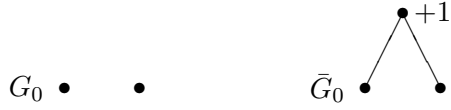
\begin{center}
\vspace*{-4ex}
\unitlength 1mm
\begin{picture}(100,30)(0,0)
\put(20,10){\circle*{1.5}}
\put(30,10){\circle*{1.5}}

\put(60,10){\circle*{1.5}}\put(70,10){\circle*{1.5}}
\put(65,20){\circle*{1.5}}
\put(60,10){\line(1,2){5}}
\put(70,10){\line(-1,2){5}}

{
\put(15,10){\makebox(0,0){$G_0$}}
\put(55,10){\makebox(0,0){$\bar G_0$}}
\put(69,20){\makebox(0,0){$+1$}}
}
\end{picture}
\vspace*{-4ex}
\end{center}

\paragraph{Brian\c{c}on polynomial.}
Let $f(x,y) = yp^3+p^2s+a_1ps +a_0s$
with $s = xy+1, \quad p = x(xy+1)+1, \quad a_1 = -\frac{5}{3}, \quad a_0=-\frac{1}{3}.$
The bifurcation set is $\B = \Binf= \{ 0, c=-\frac{16}{9} \}$, 
moreover all fibers are smooth and irreducible.
The value $0$ is not acyclic then $j_0$ is neither injective nor surjective
(but $j_\infty$ is surjective).
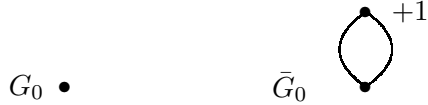
\begin{center}
\vspace*{-4ex}
\unitlength 1mm
\begin{picture}(100,30)(0,0)
\put(20,10){\circle*{1.5}}
\put(60,10){\circle*{1.5}}
\put(60,20){\circle*{1.5}}
\qbezier(60,10)(53,15)(60,20)
\qbezier(60,10)(67,15)(60,20)
{
\put(15,10){\makebox(0,0){$G_0$}}
\put(50,10){\makebox(0,0){$\bar G_0$}}
\put(66,20){\makebox(0,0){$+1$}}
}
\end{picture}
\vspace*{-4ex}
\end{center}

The value $c$ is acyclic, and $F_c$ is connected (since irreducible)
then $j_c$ is injective. The morphism $j_c$ is not surjective: $j_\infty$ is 
not surjective because the compactification of $F_c$ does not intersect the bamboo
at the last component.
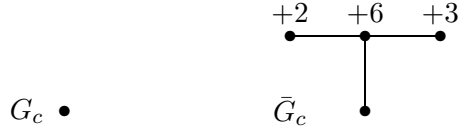
\begin{center}
\vspace*{-2ex}
\unitlength 1mm
\begin{picture}(100,30)(0,0)
\put(20,10){\circle*{1.5}}
\put(60,10){\circle*{1.5}}

\put(60,20){\circle*{1.5}}
\put(50,20){\circle*{1.5}}
\put(70,20){\circle*{1.5}}
{
\put(15,10){\makebox(0,0){$G_c$}}
\put(50,10){\makebox(0,0){$\bar G_c$}}
\put(60,23){\makebox(0,0){$+6$}}
\put(70,23){\makebox(0,0){$+3$}}
\put(50,23){\makebox(0,0){$+2$}}
\put(60,10){\line(0,1){10}}
\put(60,20){\line(1,0){10}}
\put(60,20){\line(-1,0){10}}
}
\end{picture}
\vspace*{-4ex}
\end{center}

\section{Situation around an irregular fiber}

For $f : \Cc^n \longrightarrow \Cc$ we study the neighborood of an irregular fiber.

\subsection{Smooth part of $F_c$}

Let fix a value $c \in \Cc$ and let $B^{2n}_R$ be a large closed 
ball ($R \gg 1$). Let $B^{2n}_1,\ldots,B^{2n}_p$
be small open balls around the singular points 
(which are supposed to be isolated)
of $F_c$ :  $F_c \cap \Sing$. We denote  $B^{2n}_1\cup \ldots \cup B^{2n}_p$ by $B_\cup$.
Then the \emph{smooth part} of $F_c$ is
$$\Frond_c = F_c\cap B^{2n}_R\setminus B_\cup.$$

It is possible to embed 
 $\Frond_c$ in the generic fiber  $\Fgen$ (see \cite{MW} and \cite{NN2}).
We now explain the construction of this embedding by W.~Neumann and P.~Norbury.
As $F_c$ has transversal intersection with the balls of $B_\cup$ and with
$B^{2n}_R$, then there exists a small disk $D^2_\epsilon(c)$ such that
for all $s$ in this disk, $F_s$  has transversal intersection with these balls.
According to Ehresmann fibration theorem, $f$ induces
a locally trivial fibration
$$f_|: f^{-1}\big(D^2_\epsilon(c)\big)\cap B^{2n}_R \setminus B_\cup
\longrightarrow D^2_\epsilon(c).$$
In fact, as $D^2_\epsilon(c)$ is null homotopic, this fibration is trivial.
Hence $\Frond_c \times D^2_\epsilon(c)$ is diffeomorphic to
$ f^{-1}(D^2_\epsilon(c))\cap B^{2n}_R \setminus B_\cup$.
That provides an embedding of $\Frond_c$ in  $F_s$ for all $s$ in $D^2_\epsilon(c)$;
and for such a $s$ with $s \neq c$, $F_s$ is a generic fiber.
The morphism induced in homology  by this embedding is denoted by $\ell_c$.
Let $j^\circ_c$ be the morphism induced by the inclusion of $F^\circ_c$ in 
$T_c = f^{-1}(D^2_\epsilon(c))$.
Similarly $k_c$ denotes the morphism induced by the inclusion 
of  the generic fiber $\Fgen = F_s$ (for $s\in  D^2_\epsilon(c)$, $s \neq c$) in $T_c$.
As all morphisms are induced by natural maps we have the  lemma:
\begin{lemma}
\label{prop:jkl}
The following diagram commutes:
$$
\xymatrix{
{}H_q(\Frond_c)\ar[r]^{j_c^\circ}\ar[d]_{\ell_c}    &H_q(T_c)  \\
  H_q(\Fgen)    \ar[ur]_{k_c}   & {}
}.
$$
\end{lemma}

\subsection{Invariant cycles by $h_c$}
\label{subsec:locinv}

Invariant cycles by the monodromy $h_c$ can be recovered by the following property.

\begin{proposition}
\label{prop:invcyc}
$$ \Ker\big(h_c-\id \big) = \ell_c\big(H_q(\Frond_c)\big).$$
\end{proposition}
For $n=2$, there is a similar formula in \cite{MW}, even for non-isolated singularities.

\begin{proof}
The proof uses a commutative diagram due to W.~Neumann and P.~Norbury  \cite{NN2}:
$$
\xymatrix{
{}H_q(\Fgen,\Frond_c)\ar[r]^-\sim_-\psi     &V_q(c) \ar[d]_i^{\subset} \\
  H_q(\Fgen)    \ar[u]^\varphi \ar[r]^{\id-h_c} & H_q(\Fgen)
}
$$
The morphism $i$ is the inclusion and $\psi$ is an isomorphism, so $\Ker(h_c-\id)$ 
equals $\Ker \varphi$.
The long exact sequence for the pair $(\Fgen,\Frond_c)$ 
is: 
$$ \cdots \longrightarrow H_q(\Frond_c) \stackrel{\ell_c}{\longrightarrow} H_q(\Fgen) 
\stackrel{\varphi}{\longrightarrow} H_q(\Fgen,\Frond_c) \longrightarrow \cdots$$
So $\Im \ell_c = \Ker \varphi = \Ker(h_c-\id)$.
\end{proof}

We are able to applicate this result to the calculus of the rank of $\Ker (h_c-\id)$
for $n=2$. Let denote the number of irreducible components in $F_c$ by $r(F_c)$,
and let  $\Sing_c$ be $\Sing \cap F_c$: the affine singularities on $F_c$.
Then $H_2(\Fgen,\Frond_c)$ has rank the cardinal of $\Sing_c$ which is also the
rank of  $\Ker \ell_c$. Moreover  $\rk H_1(\Frond_c)=r(F_c)-\chi(F_c)+\# \Sing_c$.
\begin{align*}
\rk \Ker \big( h_c - \id \big) 
                    &= \rk \Im \ell_c \\
                    &= \rk H_1(\Frond_c) - \rk \Ker \ell_c \\
                    &= r(F_c) - \chi(F_c) + \# \Sing_c - \# \Sing_c \\
                    &= r(F_c) - \chi(F_c).
\end{align*}

\begin{remark*}
We obtain the following fact (see \cite{MW}): 
if the fiber  $F_c$ ($c\in \B$) is irreducible then $h_c \neq \id$.
The proof is as follows: if $r(F_c)=1$ and $h_c = \id$ then from one hand
$\rk \Ker (h_c-\id) = \rk H_1(\Fgen)  = 1 - \chi(\Fgen)$ 
and from the other hand $\rk \Ker (h_c-\id) = 1 - \chi(F_c)$;
thus $\chi(F_c) = \chi(\Fgen)$ which is absurd for $c$ in $\B$ by 
Suzuki formula.
\end{remark*}

\subsection{Vanishing cycles}

Now and until the end of this paper homology is homology
with complex coefficients.

\paragraph{Vanishing cycles for eigenvalues $\lambda \neq 1$.}

Let  $E_\lambda$ be the space $E_\lambda = \Ker (h_c-\lambda\id)^p$
for a large integer $p$.

\begin{lemma}
If $\lambda \not=1$ then
$E_\lambda \subset V_q(c).$
\end{lemma}

\begin{proof}
If $\sigma \in H_q(\Fgen)$ then
$h_c(\sigma) - \sigma  \in V_q(c).$
This is just the fact that the cycle 
$h_c(\sigma) - \sigma$ corresponds to the boundary
of a ``tube'' defined by the action of the geometrical monodromy.
We remark this fact can be generalized for $j \geq 1$ to
$$h_c^j(\sigma)-\sigma \in V_q(c).$$
Let $p$ be  an integer that defines $E_\lambda$, then for $\sigma \in E_\lambda$:
\begin{align*}
 0 & = (h_c-\lambda\id)^p(\sigma)  = \sum_{j=0}^p \binom{p}{j}(-\lambda)^{p-j}h_c^j(\sigma) \\
   & = \sum_{j=0}^p \binom{p}{j}(-\lambda)^{p-j} \big(h_c^j(\sigma)-\sigma\big) 
                      + \sum_{j=0}^p \binom{p}{j}(-\lambda)^{p-j}\sigma \\
   & = \sum_{j=0}^p \binom{p}{j}(-\lambda)^{p-j} \big(h_c^j(\sigma)-\sigma\big) + (1-\lambda)^p\sigma.\\
\end{align*}
Each $h_c^j(\sigma)-\sigma$ is in $V_q(c)$, and a sum of such elements is also in
$V_q(c)$, then $ (1-\lambda)^p\sigma \in V_q(c)$. 
As $\lambda \neq 1$, then $\sigma \in V_q(c)$.
\end{proof}

\paragraph{Vanishing cycles for the eigenvalue $\lambda = 1$}
\label{subsec:valun}
We study what happens for cycles associated to the eigenvalue $1$.
Let recall that vanishing cycles $V_q(c) = \Ker k_c$ for the value $c$, are
cycles that ``disappear'' when the generic fiber tends to the fiber $F_c$.
Hence cycles that will not vanish are cycles that already exist
in  $F_c$. From the former paragraph these cycles
are associated to the eigenvalue $1$.

Let  $(\tau_1,\ldots,\tau_p)$ be a family of $H_q(\Fgen)$ such that
the matrix of $h_c$ in this family is:
$$\begin{pmatrix}
1 & 1 &  &  & (0)   \\
  & 1 & 1&  &    \\
  &   & 1 &\ddots &\\
  &(0) &  &  \ddots &1\\
  & &  &   &1\\
\end{pmatrix}.
$$

Then, the cycles $\tau_1,\ldots,\tau_{p-1}$ are vanishing cycles.
It is a simple consequence of the fact that $h_c(\sigma)-\sigma \in V_q(c)$, because for $i=1,\ldots,p-1$, we have
$h(\tau_{i+1}) - \tau_{i+1} = \tau_i$, and then
$\tau_i$ is a vanishing cycle.
It remains the study of the cycle $\tau_p$ and the particular case of Jordan blocks  $(1)$
of size $1\times 1$. 
We will start with the second part.

\paragraph{Vanishing and invariant cycles.}

Let $K_q(c)$ be  invariant and vanishing cycles for the value $c$: 
$ K_q(c)= \Ker(h_c-\id) \cap V_q(c)$.
Let us remark that the space $K_q(c)\oplus \bigoplus_{c'\not=c}{V_q(c')}$ is
not equal to $\Ker(h_c-\id)$. But  equality holds in cohomology.    

\begin{lemma}
\label{prop:inveva}
$K_q(c) = \ell_c(\Ker{j_c^\circ}).$
\end{lemma}
This lemma just follows from the description of invariant cycles 
(proposition \ref{prop:invcyc}) and from the diagram 
of lemma \ref{prop:jkl}.
For $n=2$ we can calculate the dimension of $K_1(c)$.

\begin{proposition}
\label{prop:rankkc}
For $n=2$, $\rk K_1(c) = r(F_c)-1+\rk H_1(\bar{G}_c).$
\end{proposition}

\begin{proof}
The proof will be clear after the following remarks:
\begin{enumerate}
  \item $K_1(c) = \ell_c (\Ker \jrond_c)$, by lemma \ref{prop:inveva}.
  \item $\jrond_c = j_c \circ i_c$ with $i_c : H_1(\Frond_c) \longrightarrow H_1(F_c)$
the morphism induced by inclusion. It is consequence of the commutative diagram:
$$
\xymatrix{
{} H_1(F_c) \ar[rd]^-{j_c} &  \\
 H_1(\Frond_c) \ar[u]^-{i_c} \ar[r]_-{\jrond_c} & H_1(T_c) } 
$$
  \item $\rk \Ker \jrond_c = \rk \Ker i_c + \rk \Ker j_c \cap \Im i_c$,
which is general formula for the kernel of the composition of morphisms.
  \item $\Ker j_c \cap \Im i_c = \Ker j_c$, because cycles of $H_1(F_c)$ that do not
belong to $\Im i_c$ are cycles corresponding to $H_1(G_c)$, so they already exist in
$F_c$ and are not vanishing cycles.
  \item $\rk \Ker i_c= \sum_{z\in \Sing_c} r(F_{c,z})$, where
 $F_{c,z}$ denotes the germ of the curve $F_c$ at $z$.
  \item $\rk \Ker j_c= \rk \Ker j_\infty = \sum_{P\in \Dic_c} (n(F_P)-1)
= n(F_c) + \rk H_1(\bar{G}_c)-\rk (G_c)$, it has been proved in lemma \ref{lem:jinj}.
  \item $r(F_c)  + \rk H_1(G_c) = n(F_c) + \sum_{z\in \Sing_c}( r(F_{c,z})-1)$.
This a general formula for the graph $G_c$, the number of vertices of $G_c$ is
$r(F_c)$, the number of connected components is $n(F_c)$,
the number of loops is $\rk H_1(G_c)$ and the number of edges for a
vertex that correspond to an irreducible component $F_{\mathit{irr}}$ of $F_c$ is:
$\sum_{z\in F_{\mathit{irr}}} (r(F_{\mathit{irr},z})-1)$.
  \item $\rk K_1(c)  = \rk \Ker \jrond_c - \# \Sing_c$ because
$\Ker i_c$ is a subspace of $\Ker \ell_c$ so 
$\rk K_1(c) = \rk\Ker \jrond_c - \rk \Ker \ell_c$ and the dimension of $\Ker \ell_c$ is
$\# \Sing_c$ (see paragraph \ref{subsec:locinv}).
\end{enumerate}

We complete the proof:
\begin{align*}
\rk K_1(c) 
&= \rk \ell_c (\Ker \jrond_c) & (1)  \\
&= \rk \Ker \jrond_c - \rk \Ker \ell_c &(8) \\
&=  \rk \Ker j_c \circ i_c - \# \Sing_c & (2) \text{ and } (8)  \\           
&= \rk \Ker i_c + \rk \Ker j_c \cap \Im i_c - \# \Sing_c & (3) \\
&= \rk \Ker i_c - \# \Sing_c + \rk \Ker j_c & (4) \\
&= \sum_{z\in \Sing_c} \big(r(F_{c,z}) - 1\big) +  
        n(F_c) + \rk H_1(\bar{G}_c)-\rk (G_c) & (5) \text{ and } (6)\\
&= r(F_c) - 1 + \rk H_1(\bar{G}_c). & (7) 
\end{align*}

\end{proof}

\paragraph{Filtration.}

Let  $\phi$ be the map provided by the total resolution of $f$.
The divisor $\phi^{-1}(c)$ is denoted by $\sum_i m_i D_i$ where $m_i$ 
stands for the multiplicity of $D_i$.
We associate to $D_i$ a part of the generic fiber denoted by $F_i$ (see \cite{MW}).
The \emph{filtration} of the homology of the generic fiber is the sequence of inclusions:
$$ W_{-1} \subset W_0 \subset W_1 \subset W_2=H_1(\Fgen).$$
with
\begin{itemize}
  \item $W_{-1}$: the \emph{boundary cycles}, that is to say, if
$\Fgenb$ is the compactification of $\Fgen$ and $\iota_* : H_1(\Fgen) 
\longrightarrow H_1(\Fgenb)$ is induced by inclusion then $W_{-1} = \Ker \iota_*$.
  \item $W_0$: these are \emph{gluing cycles}: the homology group
on the components of $F_i \cap F_j$ ($i\neq j$). 
  \item $W_1$: the direct sum of the $H_1(F_i)$.
  \item $W_2 = H_1(\Fgen)$.
\end{itemize}
The subspaces $W_0$ and $W_1$ depend on the value $c$.

\paragraph{Jordan blocks for $n=2$.}

For polynomials in two variables, the size of Jordan block for
the monodromy $h_c$ is less or equal to $2$.
Let denote by  $\sigma$ and $\tau$  cycles of  $H_1(\Fgen)$ such that
$h(\sigma) = \sigma$ and $h(\tau) = \sigma + \tau$.
The matrix of $h_c$ for the family $(\sigma,\tau)$ is 
$\left(\begin{smallmatrix} 1&1\\ 0&1\\
\end{smallmatrix}\right)$.
We already know that the cycle $\sigma$ vanishes.

A \emph{large cycle} is a cycle of $W_2 = H_1(\Fgen)$ that has
a non-trivial class in  $W_2/W_1$. 
According to \cite{MW} $\tau$ is large cycle;
moreover large cycles associated to the eigenvalue $1$ are
the embedding of  $H_1(\bar{G}_c)$ in $H_1(\Fgen)$.
So large cycles are not vanishing cycles.
The number of classes of large cycles in $W_2/W_1$ is 
$\rk H_1(\bar{G}_c)$, this is also the number
of Jordan $2$-blocks for the eigenvalue $1$. 

\paragraph{Vanishing cycles.}

We describe vanishing cycles.
For all the spaces  $W_{-1}$, $W_0 / W_{-1}$, $W_1/W_0$ and $W_2/W_1$
the cycles associated to eigenvalues different from $1$
are vanishing cycles.

\begin{proposition}
Vanishing cycles for the eigenvalue $1$ are dispatch as follows:
\begin{itemize}
  \item for $W_{-1}$: $r(F_c)-1$ cycles,
  \item for $W_0$: $\rk H_1(\bar{G}_c)$ other cycles,
  \item $W_1$, $W_2$: no other cycle.
\end{itemize}
\end{proposition}

Let explain this distribution. We have already remark that large cycles associated to 
$\left(
\begin{smallmatrix}  1&1\\ 0&1 \end{smallmatrix}
\right)$
are not vanishing cycles, so  vanishing cycles in 
$W_2$ are in $W_1$. Moreover there is  $\rk H_1(\bar{G}_c)$
Jordan $2$-blocks for the eigenvalue $1$ that provide
$\rk H_1(\bar{G}_c)$ vanishing cycles (like $\sigma$) in $W_0$.
The other vanishing cycles for the eigenvalue $1$ are
invariant cycles by $h_c$, in other words they belong to $K_1(c)$.
We have $W_1 \cap K_1(c) = W_0 \cap K_1(c)$ because
invariant cycles for $W_1$ that are not in $W_0$ correspond
to the genus of the smooth part  $\Frond_c$ of $F_c$
(this is due to the equality $\Ker(h_c-\id) = \ell_c(H_1(\Frond_c))$).
As they already appear in $F_c$, theses cycles are not vanishing cycles for the value $c$.
Finally, if we have two distinct cycles $\sigma$ and $\sigma'$ in $W_0 \cap K_1(c)$, with
the same class in $W_0/W_{-1}$, then $\sigma' = \sigma + \pi$, $\pi \in W_{-1}$;
this implies that $\pi = \sigma'-\sigma$ is a vanishing cycle of $K_1(c)$. 
We can choose the $r(F_c)-1$ remaining cycles of $K_1(c)$ in $W_{-1}$.


\bigskip
\bigskip
\vfill

{\noindent
Arnaud \textsc{Bodin} \\
Universit\'e Paul Sabatier Toulouse III, laboratoire \'Emile Picard, \\
118 route de Narbonne, 31062 Toulouse cedex 4, France. \\
\texttt{bodin@picard.ups-tlse.fr}
}

\end{document}